\documentclass[11pt]{amsart}
\calclayout
\usepackage{url}

\DeclareMathOperator{\End}{End}
\DeclareMathOperator{\Hom}{Hom}
\DeclareMathOperator{\tr}{tr}

\DeclareMathOperator{\coker}{coker}

\DeclareMathOperator{\ind}{ind}

\newcommand{\Spin}{\mbox{Spin}}

\DeclareMathOperator{\ch}{ch}
\DeclareMathOperator{\SO}{SO}

\DeclareMathOperator{\Str}{Str}

\newcommand{\Z}{{\mathbb Z}}
\newcommand{\ZZ}{\mathbb Z}
\newcommand{\RR}{\mathbb R}

\newcommand{\Ga}{\Gamma}

\newcommand{\hX}{\widehat{X}}
\newcommand{\hA}{\widehat{A}}
\newcommand{\hG}{\widehat{G}}
\newcommand{\hGc}{\widehat{G}_c}
\newcommand{\hP}{\widehat{P}}

\newcommand{\hg}{\widehat{g}}

\newcommand{\del}{\nabla}

\newcommand{\ci}{C^\infty}
\newcommand{\Om}{\Omega}

\theoremstyle{plain}
\newtheorem{theorem}{Theorem}[section]

\newtheorem{proposition}[theorem]{Proposition}

\theoremstyle{definition}
\newtheorem{definition}[theorem]{Definition}

\theoremstyle{remark}

\numberwithin{equation}{section}
\numberwithin{figure}{section}

\newcommand{\cE}{{\mathcal E}}

\newcommand{\cF}{{\mathcal F}}

\renewcommand{\a}{\alpha}
\renewcommand{\b}{\beta}
\renewcommand{\c}{\gamma}

\begin{document}

\title[Gerbes and the index theorem]{Gerbes, Clifford modules and the index theorem}
\author{Michael K. Murray}
\address[Michael K. Murray]
{Department of Pure Mathematics\\
University of Adelaide\\
Adelaide, SA 5005 \\
Australia}
\email[Michael K. Murray]{mmurray@maths.adelaide.edu.au}
\thanks{The author acknowledges the support of the Australian
Research Council.}
\author{Michael A. Singer }
\address[Michael A. Singer]
{School of Mathematics \\
King's Buildings
\\
University of Edinburgh EH9 3JZ \\
U.K.}
\email[Michael A. Singer]{michael@maths.ed.ac.uk}
\thanks{}


\begin{abstract}
The use of bundle gerbes and bundle gerbe modules is 
considered as a replacement for the usual theory of 
Clifford modules  on manifolds that fail to be spin. 
It is shown that both sides of the Atiyah-Singer index 
formula for coupled Dirac operators can be given natural 
interpretations using this language and that the resulting
formula is still an identity.
\end{abstract}

\maketitle

\section{Introduction}
If $M$ is an even-dimensional,  oriented, spin-manifold with
spin-bundles $S^{\pm}$ and
$W$ is a vector bundle with unitary connection $A$ then the
{\em coupled Dirac operator} is a natural first-order differential operator
\begin{equation}
\label{eq:dirac}
D_A^+ \colon \ci(M, E^+) \to \ci(M, E^-).
\end{equation}
where
\begin{equation} \label{1.27.9}
E^+:=  S^+\otimes W,\;\;  E^-:= S^-\otimes W.
\end{equation}
If $M$ is compact, the Atiyah-Singer index theorem \cite{AS0, AS1} 
gives a topological
formula for the index of $D^+_A$,
\begin{equation}
\label{eq:indextheorem}
\ind(D^+_A):= \dim \ker (D^+_A)  - \dim \coker (D^+_A) = 
\langle \hA (M) \ch(W),  [M] \rangle.
\end{equation}

If $M$ is not a spin-manifold, `generalized Dirac operators' can be
introduced, even though the Dirac operator itself is not
well-defined. This can be done by observing that $E = E^+\oplus E^-$
is a {\em Clifford module} (Clifford multiplication is extended to act
as the identity on $W$) so that a compatible connection on $E$ leads
to a Dirac operator as before.  If $M$ is not spin, there is an index
theorem for Dirac operators defined on (hermitian) Clifford modules $E$
which can be expressed in the form \cite{BerGetVer}
\begin{equation} \label{eq:2.27.8}
\ind(D^+_A) = \langle \widehat{A}(M)\ch(E/S), [M]\rangle
\end{equation}
where the notation will be explained below.  The definition is such
that if $M$ is spin and \eqref{1.27.9} holds, then the {\em relative
Chern character} $\ch(E/S)$ reduces to $\ch(W)$, the Chern character
of $W$.

In this note we explore an alternative approach to this index theorem,
using {\em bundle gerbes} and {\em bundle gerbe modules}
\cite{Mur, BouCarMatMurSte}. The point is
that there is a so-called spin bundle gerbe $\Ga$, defined on any 
manifold, which is
trivial if $M$ is spin, and modules for $\Ga$ are defined for any
representation of the spin group. In particular, the {\em spin}
$\Ga^1$-module is  always defined, and this reduces, in a suitable sense, to
the spin-bundle, if $M$ happens to be spin.  We construct a
`twisted Dirac operator' $\overline{D}$ between the two spin
$\Ga^1$-modules, and show that when this is coupled to any suitable
auxiliary $\Ga^{-1}$-module $W$ then the resulting differential 
operator reduces
to a Dirac-type operator between Clifford modules.

Furthermore, we show how to define the Chern character $\ch(W)$ of the
$\Ga^{-1}$-module $W$ and that the index formula \eqref{eq:indextheorem}
holds with this definition.  Once all the definitions are in place,
the proof of the  index theorem is mainly a question of matching up
the formalism of Clifford modules with that of
$\Ga^d$-modules.

Some of our constructions are relevant, and familiar, in relation to 
anomalies in quantum field theory. The failure of $M$ to be 
spin can be regarded as an anomaly in the global definition of spinor fields.
As is  is typically done in anomaly theory our solution is to 
introduce some additional
fields which also have an anomaly in their global definition and 
choose these so that
the two anomalies cancel in the coupled theory.

We note two recent preprints. The first
\cite{MatMelSin} discusses the
interesting situation of the families index theorem for
projective families of pseudo-differential operators. In that case the
index takes value in the twisted $K$-theory of the parameter space.
Although this is related in spirit to this note, the situation considered there
is rather different. The second \cite{Ari} presents a general Atiyah-Singer
index theorem for gerbes of which our result is an example. 

The first author thanks Mathai Varghese for useful discussions.

\section{Background on lifting bundle gerbes}

In this section we recall the definition of the lifting bundle gerbe
and  the corresponding bundle gerbe modules. The references
for this material are \cite{Mur, BouCarMatMurSte}.

\subsection{Generalities on central extensions of Lie groups}
Let
\begin{equation} \label{1.25.9}
0 \to \ZZ_k \to \hG \stackrel{\pi}{\to} G \to 1
\end{equation}
be a central extension of finite-dimensional Lie groups, where $\ZZ_k$ is
the cyclic  subgroup $\ZZ/k\ZZ$  of the circle $U(1)$. The example of 
ultimate interest
in this note is
\begin{equation} \label{4.27.9}
0 \to \ZZ_2 \to \Spin(n) \stackrel{\pi}{\to} \SO(n) \to 1,
\end{equation}
where $\Spin(n)$ is the universal cover of $\SO(n)$, but it is clearer
to work with the general case to begin with.

We will be interested in spaces on which $\hG$ and $G$ act.
Notice that if $G$ acts on a space then so also does $\hG$ via the projection.
We shall often assume  this action in the discussion below without 
further comment. In such a 
case the center $\ZZ_k$ of course acts trivially.   
To begin our discussion let $X$ and $\hX$ be right principal sets for
$G$ and $\hG$ respectively. As sets, $X$ and $\hX$ are just $G$
and $\hG$: in particular we can choose  a map $\pi:\hX \to X$ which
is compatible with the respective actions of $\hG$ on $\hX$ and
of $G$ on $X$.  Notice that  $\pi \colon \hG \to G$ induces an
action of $\hG$ on $X$.

Let $\rho:\hG\to GL(V)$ denote a (finite-dimensional) representation of $\hG$
and let $V_X$ denote the bundle $X\times V$.  With the 
$\hG$ action $\hat g(x, v) = (x\hat g^{-1}, \rho(g)v)$
this is an equivariant $\hG$ bundle.  As the isotropy 
subgroup of the $\hG$ action at any $x \in X$ is $\ZZ_k$ we must
have that $\ZZ_k$ acts on the fibres of $V_X \to X$ and it is clear
that it does so by the restriction of the representation $\rho$.

The product $X^2:= X\times X$ has a canonical map
$$
\tau :X^2 \to G,\;\;x_1\tau(x_1,x_2) = x_2,
$$
so that $\tau(x_1,x_2)$ is the unique element of $G$ that translates 
$x_1$ to $x_2$. Let $\Ga=
\tau^{*}\hG \to X^2$ be the pull-back of the $\ZZ_k$ bundle $\hG \to 
G$.   By definition
$$
\Ga_{(x_1,x_2)} = \{\hat{g}\in \hG: \pi(g) = \tau(x_1,x_2)\}
$$
is the set of all lifts to $\hG$ of $\tau(x_1,x_2)$.

The $\ZZ_k$-principal bundle $\Ga $ contains almost the same information as
$\hX$: if we denote by $\hX_1$ and   $\hX_2$
the pull-backs of $\hX\to X$ by
the two projections $X^2 \to X$, then the principal bundle of
$\ZZ_k$-equivariant maps $\hX_1 \to \hX_2$ is canonically isomorphic
to $\Ga$.  Indeed, at $(x_1,x_2)$,
the fibres of $\hX_1$ and $\hX_2$ are by
definition
$$
\{\xi: \pi(\xi) = x_1\}\mbox{ and }
\{\xi: \pi(\xi) = x_2\}
$$
respectively and it is clear that equivariant maps between these sets 
correspond
precisely to elements of $\Ga_{x_1,x_2}$.

The $\hG$-equivariant bundle $V_X\to X$ is a
`module' for $\Ga$ in the following sense. Let
$V_1$ and $V_2$ be the pull-backs of $V_X$ by the two projections
$X^2 \to X$.  Then $\Ga_{(x_1,x_2)}\subset \hG$ acts through
$\rho$ to map the fibre at $(x_1,x_2)$ of $V_1$ to that of $V_2$.

Notice that if $\rho(\ZZ_k) =\{1\}$ then $\rho$ is really a
representation of $G$ and one can divide $V_X$ by $G$ to get back
the representation space $V$. In particular,
$\End(V_X)$ is always $G$-equivariant so that $\End(V_X)/G = \End(V)$.

Finally we notice that, because $\ZZ_k$ is finite,  there is a natural $\hat
G$-equivariant flat connection $\del_X$ on any of the bundles $V_X$.

\subsection{The lifting bundle gerbe and associated modules}

Let $M$ be a smooth manifold of dimension $n$. We now globalize the
considerations of the previous section so that they correspond to the
case  when $M$ is a point.

We stay with the central extension \eqref{1.25.9} and replace $X$ by
a given principal $G$-bundle $\pi \colon P\to M$. Correspondingly 
$X^2$ is
replaced by the fibre product $P^{[2]} \to M$.  Extending the notation
from the previous section, the map $\tau:
P^{[2]}\to  G$ is defined exactly as before and we put
$\Ga = \tau^*\hG$. This is a $\ZZ_k$ bundle over $P^{[2]}$ and
in the language of \cite{Mur} is called the {\em lifting bundle gerbe}
associated to $P$ and the central extension \eqref{1.25.9}.  Note 
that $\Ga$  is a $\ZZ_k$ bundle
gerbe whereas \cite{Mur} was mostly concerned with $U(1)$ bundle gerbes.

We will be interested in $\hG$ equivariant  bundles over $P$. 
   Recall first that a finite-dimensional
representation $U$ of $\ZZ_k$ decomposes into a direct sum
of weight spaces $U_d$ defined to be the largest subspaces
on which the action of every $z \in \ZZ_k$ is just multiplication
by $z^d$. We say that $\ZZ_k$ acts on $U_d$ with weight $d$.
Notice that if $m$ is any integer and $z \in \ZZ_k$ then $z^{d+mk} = z^d$
so that the weight $d$ is an element of $\ZZ_k$.

Let   $W \to P $ be  a $\hG$ equivariant bundle for the action of $\hG$
on $P$.
  The isotropy subgroup $\ZZ_k$  acts
on the fibres of $W \to P$. If it acts with weight $d$ then
  we call $W$ a $\Ga^d$-module.
In the language of \cite{BouCarMatMurSte} a $\Ga^d$-module is precisely a
bundle gerbe module for the bundle gerbe $\Ga^d = \Ga^{\otimes d}$.

We note, as an aside, that if $W$ is any
$\hG$-equivariant bundle then at every point $p \in P$  we can decompose the
fibre of $W$ into weight spaces  for the $\ZZ_k$ action. The action
of $\Ga$ maps fibres to fibres and, because $\ZZ_k$ is central,
  preserves these weight spaces.
Hence  $W$ is a vector bundle direct sum of $\Ga^d$-modules so we lose
nothing by concentrating attention on $\Ga^d$-modules.

  Of particular importance for us are $\Ga^0$-modules. Then every 
element of $\ZZ_k$ acts by the
identity on $W$ so that  $W$ is $G$-equivariant and hence  it
descends to a bundle on $M$ after quotienting by $G$.

If  $\rho:\hG\to GL(V)$ is a representation then
$V_P = P \times V$   defines a $\hG$-equivariant vector bundle, with 
$\hG$ action $g(p, v) = (pg^{-1}, gv)$ and 
 such that the restriction of $V_P$ to each fibre of $P$ is
isomorphic to $V_X\to X$.  If the action $\ZZ_k$ defined by $\rho$ is 
of weight $d$ then $V_P$ is a $\Ga^d$-module.
Note that $V_P$ is a trivial module if and only if the
representation of $\hG$ factors through a representation $\rho'$, say,
of $G$. Then
$$
V_P = \pi^{-1}(V_M),\mbox{ where }V_M:= P\times_G V
$$
where $G$ acts by $\rho'$ on $V$ (so that $V_M$ is the vector bundle over $M$
associated by $\rho'$ to $P$).

Define a group $\hGc$ by
$$
\hGc = \left( U(1) \times \hG \right) / \ZZ_k
$$
where $\ZZ_k \subset U(1) \times \hG$ is included as the
anti-diagonal subgroup $\{(z, z^{-1}) \mid z \in \ZZ_k \}$.
The group $\hGc$ is a central extension of
$G$ by $U(1)$ and there is a short exact sequence
$$
0 \to U(1) \to \hGc \to G \to 0
$$
induced by \ref{1.25.9}.  A $\Ga^d$-module $W$ has a natural
  $\hGc$ bundle action defined by  $(z, g) w =
z^dgw$. If  it happens that there exists a principal
$\hGc$-bundle $\hP_c \to M$ covering $P$, then one can pull back any
$\Ga$-module $W$ to $\hP$ and divide by $\hG$ to get the associated
vector bundle $W_M:=\hP_c \times_{\hG} W$.  Notice that
we can repeat the construction of $\Ga$ to define a
$U(1) $-bundle $\Ga_c = \tau^{-1}(\hGc) \to P^{[2]}$.
  In \cite{Mur} it is
shown that the lift of $P$ to the $\hG_c$ bundle $\hP_c \to M$ exists 
exactly when the
$U(1)$-bundle gerbe $\hGc$ is  trivial.

Notice that there are thus two ways in which a $\Ga^d$-module $W$ may be a
legitimate vector bundle on $M$.   Either $\Ga_c$ is a trivial bundle gerbe or
$W$ is a $\Ga^0$-module, i.e $d=0$ modulo $k$.

In the next few sections
we shall make definitions (for example of connections and curvature)
for $\Ga^d$-modules on $M$ that run parallel to the usual theory of vector
bundles on $M$. The reader should have no difficulty in verifying that
these definitions always reduce to the standard ones in the case
that $d=0$ is trivial or $\Ga_c$ is trivial.

If $W$ is a $\Ga$-module then
  $\End(W)$ is a  trivial module. This situation
will occur frequently below, and we shall write
$$
\cE(W) = \End(W)/G.
$$

\subsection{Concrete representation}
We can obtain a \v{C}ech description of the above structures by choosing
a good cover
  $\{ U_\a \}$ of $M$ and
local cross-sections $s_\a \colon U_\a \to P$. Then the transition
functions
$$
g_{\a\b} \colon U_\a \cap U_\b \to G
$$
are defined by the conditions $s_\a = s_\b g_{\a\b}$. The problem of
constructing a covering principal bundle $\hat{P}$ is the same as
that of making a consistent choice of lifts $\hg_{\a\b}$ of the
$g_{\a\b}$. Indeed, if we make an arbitrary choice of such lifts, then
we obtain a $\ZZ_k$-valued \v{C}ech 2-cocycle $ e_{\a\b\c} \in H^2(M, 
\ZZ_k)$ defined as
\begin{equation} \label{7.27.9}
  e_{\a\b\c}:=  \hat g_{\a\b}\hat g_{\b\c} \hat g_{\c\a}
\end{equation}
We can also regard $e_{\a\b\c}$ as a $U(1)$
2-cocycle by inclusion of $\ZZ_k$ in $U(1)$. Using the exact
sequence
$$
0 \to \ZZ \to \RR \to U(1)
$$
we obtain, in the usual way a class $[e] \in H^3(M, \ZZ) \simeq H^2(M, U(1))$.
This class is   the obstruction to the
existence of $\hP_c$.

Similarly, given a $\Ga^d$-module $W$,
set $W_\a = s_\a^{-1}(W)$.
The lifted transition map $\hat g_{\a\b}$ induces a
map $\phi_{\a\b} \colon W_\b \to W_\a$.
  The $\phi$'s determine a cocycle
$$
\phi_{\a\b}\phi_{\b\c} \phi_{\c\a} = e^d_{\a\b\c}
$$
If this cocycle is a coboundary then the locally defined bundles
$W_\a$ patch together to form a genuine bundle on $M$
whose pull-back is $W$.
This is  the case   if  $(e_{\a\b\c})$
is  a coboundary (as a $U(1)$ 2-cocycle) or if $d=0$.

An alternative description of the local structure of $W$ can
be obtained by defining $U(1)$ bundles
$L_{\a\b} = g_{\a\b}^{-1}(\hG_c) \to U_\a \cap U_\b$. Then there are 
isomorphisms
$W_\a = L^d_{\a\b}\otimes W_\b$.

\subsection{Operations on $\Ga^d$-modules}

The following observations are easy but important \cite{BouCarMatMurSte}:

\begin{enumerate}\item There is a bijection between
$\Ga^0$-modules and bundles on $M$.  (Induced by pull back to $P$ in one
direction, division by $G$ in the other.)

\item The tensor product of a $\Ga^{d_1}$-module and a
$\Ga^{d_2}$-module is a $\Ga^{d_1+d_2}$-module.

\item The endomorphism bundle of an $\Ga^d$-module is a  $\Ga^0$-module.

\item The direct sum of two $\Ga^d$-modules is a $\Ga^d$-module.
\end{enumerate}

These properties (especially the last one) lead to a definition of
twisted $K$-theory, by introducing the $K$-group of the semi-group of
$\Ga$-modules \cite{BouCarMatMurSte}. In this case the `twist' is by
$d.\delta[e] \in H^3(M,\Z)$ where  $[e]$ was defined in
\eqref{7.27.9}.

\subsection{The spin-bundle gerbe}
\label{10.27.9}
The main example in this paper arises by taking
the central extension
$$
0 \to \ZZ_2 \to  \Spin(n) \to SO(n) \to 1,
$$
and $P$ to be the bundle of oriented orthonormal frames in $TM$.  The
lifting bundle gerbe $\Ga$ in this case is called the {\em spin-bundle gerbe}.

In the spin case the construction of $\hGc$ is the well-known 
construction of the
group $\Spin_c(n)$. The manifold $M$ is spin when  $\Ga$ is a trivial
$\ZZ_2$ bundle gerbe  and  spin-c when $\Ga_c$ is
a trivial $U(1) $  bundle gerbe.

If (as we assume throughout) $n$ is even then we have the two
half-spin representations $\Sigma^\pm$ of $\Spin(n)$ as well as the
full spin representation $\Sigma = \Sigma^+\oplus \Sigma^-$. The
$\Ga^1$-modules associated to these representations are called the {\em
spin} $\Ga^1$-modules, and will be our replacement for honest
spin-bundles on $M$.

\section{Connections and curvature for $\Ga^d$-modules}

\subsection{Definition of $\Ga^d$-module connection}

In this note we shall generally be identifying connections with the
corresponding covariant derivative operators. In particular a
connection on $W$ is a linear differential operator
$$
\del_A:\Om^0(P,W) \to \Om^1(P,W)
$$
which satisfies the Leibnitz rule
$\del_A(fs) = df\otimes s +  f\otimes \del_A s$
for all functions $f$ and sections $s$.

\begin{definition} The connection $\del_A$ on a $\Ga$-module is called a {\em
$\Ga$-module  connection} if it is $\hGc $-equivariant and
restricts to the flat connection $\del_X$ on the fibres of $P$.
\end{definition}

If $U$ is a sufficiently small open set of $M$ then we can identify
the restriction of $P$ to $U$ with $U\times X$ and the restriction
of $W$ with $W_X$ (or rather with the pull-back of this bundle by
the projection $U\times X\to X$). With these identifications, a
$\Ga$-module connection takes the form
$$
\nabla_A = \nabla_X + \nabla_B
$$
where $\nabla_B$ is the pull-back of a connection from the base $U$.
It is straightforward to use a partition of unity argument  to patch
such local descriptions to show that $\Ga$-module connections
exist on any $\Ga$-module. The following proposition gathers some basic
information about $\Ga$-module connections.

\begin{proposition}
Let the notation be as above and $W$ be a $\Ga^d$-module.

\begin{enumerate}
\item The set of all $\Ga^d$-module connections on $W$ is an
affine space relative to the space $\Om^1(M,\cE(W))$.

\item The curvature
of a $\Ga^d$-module connection descends to define an element
$\cF_A\in\Om^2(M,\cE(W))$.

\item If $W$ and $W'$ are two $\Ga^d$-modules, and $\del_A$ and
$\del_{A'}$ are $\Ga^d$-module conections on $W$ and $W'$, then the
curvature of the tensor product connection $\del_B$ on $W\otimes W'$ is
$$
\cF_B = \cF_A \otimes 1 + 1 \otimes \cF_{A'}  $$
an element of $\Om^2(M,\cE(W\otimes
W')) = \Om^2(M,\cE(W)\otimes \cE(W'))$.

\end{enumerate}
\end{proposition}

\begin{proof} Let $\del_1$ and $\del_2$ be two connections on $W$.
It follows as usual from the Leibnitz rule that $a=\del_2-\del_1$ lies
in the vector space $\Om^1(P,\End(W))$. If now $\del_1$ and
$\del_2$ are $\Ga^d$-module connections it follows that $a$
is $\hG$-equivariant. Since the centre acts trivially here, $a$ is
actually $G$-equivariant. But from the local description it follows
that $i_\xi(a) = 0$ for any vertical vector field $\xi$. Hence we can
divide by $G$ to realize $a$ as an element of $\Om^1(M,\cE(W))$.   It
is even easier to see that if $\del_A$ is a $\Ga$-module connection
and $a\in \Om^1(M,\cE(W))$, then $\del_A + a$ is another $\Ga$-module
connection.

The proof of the statement about the curvature is very similar. The
curvature $F_A$ of $\nabla_A$ is an element of the space
$\Om^2(P,\End(W))$. This element is $G$-equivariant and satisfies
$i_\xi(F_A) = 0$ if $\xi$ is any vertical tangent vector (from the
local description and the fact that $\del_X$ is flat). Hence it
descends to an element of $\Om^2(M,\cE(W_P))$ as required.

The proof of the last statement is trivial.

\end{proof}

If $d=0$ and $W = \pi^{-1}(W_M)$ and $W_M = W/G$ then
it is easy to see that there is a bijective correspondence between
connections on $W_M$ and $\Ga^0$-module connections on $W$. Similarly
the descended curvature in (ii) above is just the curvature of the
connection on $W_M$.

We remark that if a $\Ga^d$-module $W$ is given any $\hG$ invariant 
connection $\nabla$
we can find a $\Ga^d$-module connection by choosing a connection $A$ 
on $P \to M$.
Indeed  $A$ is a one-form on $P$ with values in the Lie algebra of $G$ which
equals the Lie algebra of $\hG$ which  acts on the fibres
of $W$. The connection $\nabla - A$ is a $\Ga^d$-module connection. 
In particular
if $V$ is a representation of $\hG$ and we take the flat connection 
$d$ on $W = P \times V$  then
$d - \rho(A)$ is a $\Ga^d$-module connection. Here we abuse notation: 
$\rho$ is a
representation of $\hG$ on $V$ and we use this to obtain a 
representation of the
Lie algebra of $\hG$ and hence the Lie algebra of $G$ on $V$.

\subsection{(Twisted) Chern character}

Let the data be as before, so that $\del_A$ is a $\Ga$-module connection
on $W$, with (descended) curvature-form $\cF_A$.

\subsubsection{Definition} The (twisted) Chern character of $\cF_A$ is
defined to be
$$
\ch(\cF_A)  = \tr\exp(\cF_A/2\pi i).
$$

\vspace{10pt}
As in \cite{BouCarMatMurSte} one shows that $\ch(\cF_A/2\pi i)$
is closed and that its de Rham cohomology class is independent of the
choice of $\Ga^d$-module connection $A$.  Define $\ch(W)$ to be this
cohomology class, called the (twisted) Chern character of the
$\Ga^d$-module $W$.  It is readily checked that if $W = \pi^{-1}(W_M)$ is a
trivial $\Ga^d$-module then $\ch(W) = \ch(W_M)$.

\subsection{(Twisted) Dirac operator}

We apply the above theory in our main example, to construct a
replacement for the Dirac operator. Recall the notation of
\S\ref{10.27.9} for the spin-bundle gerbe, and denote by
$\RR^n$ the fundamental representation of $\SO(n)$, so that
$\RR^n_P  = \pi^*(TM)$ is the pull-back to the frame-bundle of the tangent
bundle. Then Clifford multiplication is defined as a
map between $\Ga$-modules
$$
c:(\RR^n_P)^*\otimes \Sigma_P^+ \to \Sigma_P^-
$$
The Riemannian metric on $M$ induces a standard connection on the
frame-bundle $P$ and this in turn induces $\Ga$-module connections on
all associated bundles.
Using these data, we define a differential operator (the twisted Dirac
operator)
\begin{equation} \label{1.26.9}
\overline{D}^+:\ci(P,\Sigma^+_P) \to \ci(P,\Sigma^-_P)
\end{equation}
by composing covariant differentiation with Clifford
multiplication. This operator is  $\Spin(n)$-equivariant but not
elliptic  because it doesn't differentiate in the vertical directions.

Now let $W$ be a $\Ga^{-1}$-module with
$\Ga^{-1}$-module connection $A$\footnote{Of course, as $d=2$, we have $1=-1$.}.
  Then there is a coupled version of \eqref{1.26.9}
\begin{equation} \label{2.26.9}
\overline{D}_A \colon \ci(P, \Sigma_P^+\otimes W) \to \ci(P, 
\Sigma_P^-\otimes W).
\end{equation}
(Here Clifford multiplication is extended to act trivially in $W$.)
Because $\Sigma_P^{\pm}\otimes W$ is a $\Ga^0$-module it descends
to a  bundle $E^{\pm}$ on $M$ and so also does the Dirac operator to give the
twisted Dirac operator:
\begin{equation} \label{3.26.9}
{D}^+_A
\colon \ci(M, E^+) \to \ci(M, E^-).
\end{equation}
If $M$ is spin, then both \eqref{1.26.9} and \eqref{2.26.9} descend to
$M$ and $E^\pm = \Sigma^\pm_M\otimes W_M$. In general, however, this
tensor product decomposition exists only locally on $M$.

We can now state:

\begin{proposition}[Twisted Index Theorem]
\label{prop:twisted}
If $M$ is compact, then the
Dirac operator defined above satisfies the index formula
$$
\ind(D^+_A) = \langle \hA (M) \ch(W) , [M]\rangle .
$$
\end{proposition}
The proof will be given in the next section, by showing how the
present gerbe framework fits with the more standard formalism of
Clifford modules.

\section{Clifford modules and generalized Dirac operators}

In the present section we recall the definition of Clifford modules,
and show how they correspond to $\Ga$-modules (from now on, $\Ga$ is
the spin-bundle gerbe). The index theorem for generalized Dirac
operators, in the form stated in  \cite{BerGetVer} then yields
Proposition~\ref{prop:twisted} at once.

\subsection{Clifford algebras and their representations}

If $V$ is a real vector space of even dimension $n$, equipped with a
positive-definite inner product, then $C(V)$ will denote the
corresponding (complex) Clifford algebra. This is a {\em superalgebra}
in the sense that the underlying vector space is $\ZZ_2$-graded,
$$
C(V) = C_+(V)\oplus C_-(V)
$$
  and
the algebra operations are compatible with this grading (even times
even and odd times odd are even, even times odd is odd).  If $C(V)$ is
identified with the complex exterior algebra of $V$ (but with a
different product) then $C_+(V)$ corresponds to the forms of even degree
and $C_-(V)$ corresponds to the forms of odd degree. The {\em
chirality operator} will be denoted by $\gamma$, so $\gamma = \pm 1$
on $C_\pm(V)$.  Recall that there is an embedding $c:V\subset C(V)$
called {\em Clifford multiplication}.

A $\ZZ_2$-graded complex inner product space $E$ is called a hermitian
Clifford module if there is an even action of $C(V)$ on $E$ (i.e. the
grading is preserved by $C^+(V)$ and altered by $C^-(V)$) and $c(v)$
is skew-adjoint for
all $v\in V$. The spin representation $\Sigma =
\Sigma^+\oplus \Sigma^-$ is a hermitian Clifford module which has the
property
$$
C(V) = \End(\Sigma).
$$
Moreover, it is a basic fact that every finite-dimensional hermitian 
Clifford module $E$
has the form
\begin{equation} \label{12.27.9}
E= \Sigma\otimes W
\end{equation}
where $W$ is some complex vector space carrying a trivial
$C(V)$-action. In this situation one has the basic formulae:
$$
W = \Hom_{C(V)}(S,E),\;\;\End(W) =\End_{C(V)}(E).
$$
The subscript $C(V)$ here means `commutes with the action of $C(V)$'.

\subsection{Clifford bundles and $\Ga$-modules}

Let $M$ be an even-dimensional oriented Riemannian manifold.
Introduce $C(M)$, the bundle of complex Clifford algebras of $T^*M$ and
continue to denote Clifford multiplication by $c:T^*M \to C(M)$.

A complex $\Z_2$-graded hermitian vector bundle $E= E_+\oplus E_-$  over $M$ is
called a {\em Clifford module} if
\begin{enumerate}
\item the sub-bundles $E_+$ and $E_-$ are orthogonal with respect to
       the inner product;
\item for each point $x\in M$, $E_x$ is a hermitian Clifford module
       for $C_x(M)$.
\end{enumerate}

Let $\Ga$ be the spin-bundle gerbe of \S\ref{10.27.9} and recall that
$\pi:P\to M$ is the frame-bundle. Let $E_P = \pi^{-1}(E)$ be the pull-back
of $E$ to $P$. Each point $p$  of $P$ is an isomorphism
of $\RR^n$ with the $T_{\pi(p)}M$ which maps the standard inner
product to the metric on the tangent space. It follows that $p$ induces an
isomorphism  $(E_P)_p$ can be made into a $C(\RR^n)$ Clifford module.
We define
$$
W_p = \Hom_{C(\RR^n)}(\Sigma, (E_P)_p)
$$
and we have
$$
(E_P)_p = \Sigma \otimes W_p
$$
or globally
$$
E_P = \Sigma_P \otimes W.
$$
The action of $\Spin(n)$ on $\Sigma$ induces an action
on $W$ making it a $\Ga^{-1}$-module. From the definition of $W$ we see that
it is a $\Ga^{-1}$-module. This gives
us a generalisation of the basic result of \cite{BerGetVer}:

\begin{proposition}
Let  $M$ be a Riemannian manifold with $\Ga$ its
Spin bundle gerbe. If $\Sigma_P$ is the spin bundle gerbe module then 
every Clifford module  $E$ on $M$ has the form  $E=\Sigma_P \otimes 
W$ for some $\Ga^{-1}$-module $W$.
\end{proposition}

\vspace{10pt}
\noindent{\bf Remark} In \cite{BerGetVer} it was shown that if $M$ is spin then
every Clifford module is the tensor product of the spin bundle with 
an arbitrary
bundle.

\subsection{Generalized Dirac operators}

Given a Clifford module $E$ for $C(M)$, the next step in the definition
of a generalized Dirac operator  is the introduction of a connection
$\nabla_A$ on $E$ compatible with the grading, the  metric $h$ and with
the Clifford   action in the sense that
$$
\nabla_A[c(\xi)f] = c(\nabla\xi)e + c(\xi)\nabla_A e
$$
for all $\xi\in \Omega^1(M)$ and $e\in \ci(M,E)$. In other words, if
$c$ is viewed as a smooth section of $TM\otimes \End(E)$, then $c$ is
parallel with respect to the tensor product connection on this bundle.

Given a Clifford module with a compatible connection $\nabla_A$, the
associated Dirac operator $D_A$ is defined as the composite
$$
\ci(M,E) \stackrel{\nabla_A}{\to} \ci(M,T^*M\otimes E)
\stackrel{c}{\to} \ci(M, E).
$$
Since Clifford multiplication by $T^*M$ is in the odd part of $C(M)$,
and the connection is even, the composite $D_A$ is an odd operator.
It can therefore be broken into two pieces, traditionally written
$$
D_A^{\pm}: \ci(M, E^{\pm}) \to \ci(M, E^{\mp})
$$
with the property that $D_A^+$ and $D_A^-$ are formal adjoints of each
other with respect to the standard $L^2$ inner product on $\ci(M, E)$.

In order to write a cohomological formula for  the {\em index} of $D_A^+$
we need to know a little about the curvature $F_A$ of our connection
$\nabla_A$. Now $F_A$ is a $2$-form with values in $\End(E)$ and we
can try to decompose $F_A$ according to the decomposition
$$
\End(E) = C(M)\otimes \End_{C(M)}(E).
$$
Because of the compatibility between $\nabla_A$ and $c$, we have
$$
[F_A , c(\xi) ] = c(R(\xi))
$$
for any tangent vector $\xi$. It is easy to show that there is
an element $c(R) \in C(M)$ such that
$$
[c(R) , c(\xi)] = c(R(\xi))
$$
and hence if  we define
$$
F_{E/S} = F_A - c(R)
$$
then $F_{E/S}$ automatically commutes with the action of $C(M)$. The
quantity $F_{E/S}$ is called the {\em twisting curvature} in
\cite{BerGetVer}.

If $M$ is spin and $E =\Sigma\otimes W$, the connection on
$E$ would be given by the tensor product of the spin connection on
$\Sigma$ and an auxiliary unitary connection $A$ on $W$ and the curvature
decomposes according to the tensor product; in this case, $F_{E/S} = 
1\otimes F_A$
where $F_A$  just the curvature of the auxiliary connection $A$.

If $M$ is not spin, we pass to the corresponding $\Ga$-modules,
$\Sigma_P$ and $W$. It is easy to check that a Clifford-module
connection on $E$ always lifts as the tensor product of the standard
$\Ga$-module connection on $\Sigma_P$ with a unitary $\Ga$-module
connection on $W$. Conversely the tensor product of such
connections always descends to define a Clifford-compatible connection
on $E$. Moreover, if
$\cF$ is the (descended) curvature of $\Sigma_P\otimes W$ then
we have $\cF = R \otimes 1  + 1 \otimes \cF_A$ where $R$ is the (descended)
curvature of $\Sigma_P$ (essentially the Riemann curvature of the
given Riemannian metric). Hence $F_{E/S} = 1 \otimes \cF_A$.

The index formula for generalized Dirac operators now takes the
form \cite{BerGetVer}:
\begin{equation}
\label{bgv}
\ind(D^+_A) = \langle \widehat{A}(M)\ch(E/S), [M]\rangle
\end{equation}
where the {\em relative Chern character} is defined as the relative
supertrace of the exponential of the twisting curvature:
$$
\ch(E/S) = \Str_{E/S}[\exp(-F_{E/S}/2\pi i)].
$$
Here the relative supertrace is a way of defining the trace of an
endomorphism of the auxiliary bundle $W$, without having to separate
it off. More precisely, for an endomorphism $\phi$ of $E$
(super-)commuting with the action of $C(M)$ we define
$$
\Str_{E/S}(\phi) = 2^{-n/2}\Str_E(\gamma \phi)
$$
where $\gamma$ is the chirality operator of $E$.

Using the fact that $F_{E/S} = 1 \otimes \cF_A$ we
see that $\ch(E/S) = \ch(W)$ and hence our twisted
index theorem follows.

\section{Concluding remarks}
In the discussion of lifting bundle gerbes we specialised to the
case of a central extension with center $\ZZ_k$ as for the
application to the Dirac operator we only needed $\ZZ_2$.  If
we allow central extensions with center $U(1)$ the only additional
complication is that the bundle gerbe does not have a canonical
flat connection and we need to introduce the notion of a bundle
gerbe connection and its associated curving and three-curvature \cite{Mur,
BouCarMatMurSte} which are fed into the definition of bundle gerbe module
and its Chern character. This more general case might be useful for 
constructing
a twisted Dirac-Raymond operator on a non-string manifold.

\end{document}